\documentclass[12pt, a4paper]{article}
\usepackage{amssymb}
\usepackage{latexsym}
\usepackage{amsbsy}

\newcommand{\bcen}{\begin{center}}     \newcommand{\ecen}{\end{center}}
\newcommand{\bay}{\begin{array}}      \newcommand{\eay}{\end{array}}
\newcommand{\beq}{\begin{eqnarray*}}      \newcommand{\eeq}{\end{eqnarray*}}
\def\az{\alpha}
\def\bz{\beta}

\def\gz{\gamma}

\def\sz{\sigma}

\def\lra{\longrightarrow}

\def\Aut{\mbox{Aut}}

\def\dim{\mbox{dim}}
\def\End{\mbox{End}}
\def\Ext{\mbox{Ext}}

\def\Hom{\mbox{Hom}}
\def\Im{\mbox{Im}}

\def\la{\langle}
\def\mod{\mbox{mod}}

\def\ra{\rangle}

\def\rank{\mbox{rank}}

\def\Supp{\mbox{Supp}}

\begin{document}

\title{\bf {\Large Is tame open?\footnote{Project 10201004 supported by NSFC.}}}

\author{Yang Han}

\date{\footnotesize Institute of Systems Science, Academy of Mathematics and Systems
Science,\\ Chinese Academy of Sciences, Beijing 100080,
P.R.China.\\ E-mail: hany@iss.ac.cn}

\maketitle

{\it {\footnotesize Dedicated to Professor Claus Michael Ringel on
the occasion of his 60th birthday}}

\bigskip

\bigskip

\begin{abstract}

\medskip

Is tame open? No answer so far. One may pose the Tame-Open
Conjecture: Tame is open. But how to support it? No effective way
to date. In this note, the rank of a wild algebra is introduced.
The Wild-Rank Conjecture, which implies the Tame-Open Conjecture,
is formulated. The Wild-Rank Conjecture is improved to the
Basic-Wild-Rank Conjecture. A covering criterion on the rank of a
basic wild algebra is given, which can be effectively applied to
verify the Basic-Wild-Rank Conjecture for concrete algebras. It
makes all conjectures much reliable.

\medskip

\end{abstract}

\medskip

2000 Mathematics Subject Classification: 16G60, 16G10, 16G20

\bigskip

\bigskip

Throughout $k$ denotes a fixed algebraically closed field. By an
algebra we mean a finite-dimensional associative $k$-algebra with
identity. By a module we mean a left module of finite
$k$-dimension except in the context of covering theory. We denote
by mod$A$ the category of finite-dimensional left $A$-modules. For
terminology in the representation theory of algebras we refer to
[2] and [31].

\medskip

\bcen{\bf 1. Tame-Open Conjecture}\ecen

\medskip

For $d \in \mathbb{N}_1:=\{1,2,3,...\}$, ${\cal A}_d$ denotes the
affine variety of associative algebra structures with identity on
$k^d$ (cf. [11; \S 2.1]). The linear group $GL_d(k)$ operates on
${\cal A}_d$ by transport of structure (cf. [11; \S 2.2]). One
remarkable result in the geometry of representations is {\it
finite representation type is open}, i.e., all $d$-dimensional
$k$-algebras of finite representation type form an open subset of
${\cal A}_d$ (cf. [11, 24, 13]). Inspired by this, Geiss asked
whether tame is open (cf. [13, 14])? Of course one may pose a
conjecture as follows:

\medskip

{\bf Tame-Open Conjecture.} {\it For any $d \in \mathbb{N}_1$, all
tame algebras in ${\cal A}_d$ form an open subset of ${\cal
A}_d$.}

\medskip

How to support the Tame-Open Conjecture? An obvious way is to
verify it for each dimension $d$. In the cases of $1 \leq d \leq
3$, ${\cal A}_d=\{$all $d$-dimensional tame algebras$\}$. Thus
Tame-Open Conjecture holds for $1 \leq d \leq 3$. In the case of
$d=4$, one can easily determine the representation type of all
$4$-dimensional algebras listed in [11; \S 5]. Apply the upper
semi-continuity of the function $A \mapsto \dim_k \Aut(A)= \dim_k
\End(A)$ (cf. [24; Proposition 6.3]), one can show that Tame-Open
Conjecture holds for $d=4$ as well. However, for $d \geq 5$, even
for $d=5$ only, the problem becomes too complicated to be dealt
with (cf. [18; 28]). Thus it seems that it is difficult to go
further along this way.

Note that the Tame-Open Conjecture was also studied by Kasjan from
the viewpoint of model theory. He proved that the class of tame
algebras is axiomatizable, and finite axiomatizability of this
class is equivalent to the Tame-Open Conjecture (cf. [20]).
Nevertheless, it seems that this cannot support Tame-Open
Conjecture.

\medskip

\bcen{\bf 2. Wild-Rank Conjecture.}\ecen

\medskip

A finite dimensional $k$-algebra $A$ is called {\it wild} if there
is a finitely generated $A$-$k \langle x,y \rangle$-bimodule $M$
which is free as a right $k \langle x,y \rangle$-module and such
that the functor $M \otimes _ {k \langle x,y \rangle} -$ from
$\mod k \langle x,y \rangle$ to $\mod A$ preserves
indecomposability and isomorphism classes (cf. [6]). We say that
$A$ is {\it strictly wild} if in addition the functor $M
\bigotimes_{k \la x,y \ra}-$ is full. In a natural way, we can
define the wildness or strictly wildness for a full subcategory of
the module category over an algebra. If the algebra $A$ is wild
then we denote by $r_A$ the number $\min\{\rank _{k \langle x,y
\rangle}M|M$ is a finitely generated $A$-$k \langle x,y
\rangle$-bimodule which is free as a right $k \langle x,y
\rangle$-module and such that the functor $M \otimes _ {k \langle
x,y \rangle} -$ from $\mod k \langle x,y \rangle$ to $\mod A$
preserves indecomposability and isomorphism classes$\}$. By [5;
Corollary 2.4.3], $k \langle x,y \rangle$ is a free ideal ring. By
[5; Corollary 1.1.2], $k \langle x,y \rangle$ is an IBN ring. Thus
the rank of a free $k \langle x,y \rangle$-module is unique. Hence
$r_A$ is well-defined and called the {\it rank} of the wild
algebra $A$. Similarly we may define the rank $r_{\cal C}$ of a
wild subcategory ${\cal C}$ of mod$A$. Obviously $r_A \leq r_{\cal
C}$.

In this paper, we do not distinguish the $d$-dimensional algebras
from the points in ${\cal A}_d$. Put ${\cal T}_d:=\{A \in {\cal
A}_d|A$ tame$\}$ and ${\cal W}_d:=\{A \in {\cal A}_d|A$ wild$\}$,

\medskip

{\bf Wild-Rank Conjecture.} There is a function $f: \mathbb{N}
\rightarrow \mathbb{N}$ such that $r_A \leq f(d)$ for all $A \in
{\cal W}_d$.

\medskip

{\bf Remark 1.} In some sense, the Wild-Rank Conjecture is an
analogue of the numerical criterion of finite representation type
(cf. [3; Theorem]).

\medskip

If an algebraic group $G$ acts on a variety $X$ then the {\it
number of parameters} of $G$ on $X$ is $\dim_GX:=\max\{\dim
X_{(s)}-s|s \geq 0\}$ where $X_{(s)}$ is the union of the orbits
of dimension $s$ (cf. [19; p.71] or [25; p.125] or [7; p.399]). If
$A$ is a finite dimensional $k$-algebra then the set $\mod(A,n)$
of the $n$-dimensional representations of $A$ is the closed subset
of $\Hom_k(A,M(n,k))$ consisting of all $k$-algebra homomorphisms
from $A$ to the algebra $M(n,k)$ of $n \times n$ matrices. There
is a natural conjugation action of $GL_n(k)$ on $\mod(A,n)$. Put
${\cal A}_{d,\leq n}:=\{ A \in {\cal A}_d |
\dim_{GL_n(k)}\mod(A,n) \leq n\}$ and ${\cal A}_{d,>n}:=\{ A \in
{\cal A}_d|\dim_{GL_n(k)}\mod(A,n)
> n\}$.

\medskip

{\bf Lemma 1.} ([13; Proposition 1], [7; Proof of Theorem B]) {\it
${\cal A}_{d,\leq n}$ is an open subset of ${\cal A}_d$ and ${\cal
A}_{d,>n}$ is a closed subset of ${\cal A}_d$ for all $d$ and
$n$.}

\medskip

Put ${\cal A}_d^{\leq n}:=\cap_{i=1}^n{\cal A}_{d,\leq i}$ and
${\cal A}_d^{>n}:=\cup_{i=1}^n{\cal A}_{d,>i}$. Then ${\cal
A}_d^{\leq 1} \supseteq {\cal A}_d^{\leq 2} \supseteq \cdots$ and
${\cal A}_d^{>1} \subseteq {\cal A}_d^{>2} \subseteq \cdots$. By
Lemma 1, ${\cal A}_d^{\leq n}$ is an open subset of ${\cal A}_d$
and ${\cal A}_d^{>n}$ is a closed subset of ${\cal A}_d$ for all
$d$ and $n$.

\medskip

{\bf Lemma 2.} ([9; Proposition 2], [13; Proposition 2], [7; Lemma
3]) {\it ${\cal T}_d=\cap_{i \in \mathbb{N}_1}{\cal A}_{d,\leq i}
= \cap_{i \in \mathbb{N}_1} {\cal A}_d^{\leq i}$ and ${\cal
W}_d=\cup_{i \in \mathbb{N}_1}{\cal A}_{d,>i} = \cup_{i \in
\mathbb{N}_1} {\cal A}_d^{>i}$.}

\medskip

{\bf Theorem 1.} {\it The Wild-Rank Conjecture implies the
Tame-Open Conjecture.}

\medskip

{\bf Proof.} If the Wild-Rank Conjecture holds then there is a
function $f: \mathbb{N} \rightarrow \mathbb{N}$ such that $r_A
\leq f(d)$ for all $A \in {\cal W}_d$ and $d \in \mathbb{N}_1$.
Let $A \in {\cal W}_d$. Then there is a finitely generated $A$-$k
\langle x,y \rangle$-bimodule $M$ which is free of rank $r_A$ over
$k \langle x,y \rangle$ such that the functor $M \otimes _ {k
\langle x,y \rangle} -$ from $\mod k \langle x,y \rangle$ to $\mod
A$ preserves indecomposability and isomorphism classes. Note that
$\phi := M \otimes _{k \langle x,y \rangle}-: \mod(k \langle x,y
\rangle,t) \rightarrow \mod(A,r_At)$ is a regular map (cf. [8;
p.67]). Consider the stratifications $\mod(k \langle x,y
\rangle,t)= \cup_i \mod(k \langle x,y \rangle,t)_{(i)}$ and
$\mod(A,r_At)= \cup_j\mod(A,r_At)_{(j)}$. Since $\mod(k \langle
x,y \rangle,t)$ is irreducible and $\mod(k \langle x,y \rangle,t)=
\cup_{i,j}(\mod(k \langle x,y \rangle,t)_{(i)} \cap
\phi^{-1}(\mod(A,r_At)_{(j)}))$, there are $i$ and $j$ such that
the constructible subset $X:=\mod(k \langle x,y \rangle,t)_{(i)}
\cap \phi^{-1}(\mod(A,r_At)_{(j)})$ is irreducible and dense in
$\mod(k \langle x,y \rangle,t)$. Thus $\phi(X)$ is an irreducible
and constructible subset of $\mod(A,r_At)_{(j)}$. Consider the
restriction of $\phi$ on $X$ and $\phi(X)$. By [29; \S I.8 Theorem
3], $\dim \phi(X)- \dim X= \dim \phi^{-1}(y)$ for some $y \in
\phi(X)$. Take any $x \in \phi^{-1}(y)$. Since the inverse image
of an orbit under $\phi$ is an orbit, $\phi$ induces a regular map
$\psi$ from the orbit $GL_t(k) \cdot x$ to the orbit $GL_{r_At}(k)
\cdot y$. Apply [29; \S I.8 Theorem 3] again, we have $\dim
\phi^{-1}(y) = \dim \psi^{-1}(y) = \dim GL_{r_At}(k) \cdot y -
\dim GL_t(k) \cdot x= j - i$. Therefore $\dim_{GL_{r_At}(k)}
\mod(A,r_At) \geq \dim \mod(A,r_At)_{(j)} -j \geq \dim \phi(X) -j
= \dim X +(j-i)-j = \dim \mod(k \langle x,y \rangle,t) - i > \dim
\mod(k \langle x,y \rangle,t)-\dim GL_t(k)=2t^2-t^2=t^2$ for all
$t$. In particular, take $t=r_A$ then $\dim_{GL_{r_A^2}(k)}
\mod(A,r_A^2)>r_A^2$. This implies that for any $A \in {\cal
W}_d$, $A \in {\cal A}_{d,>r_A^2} \subseteq {\cal A}_d^{>r_A^2}
\subseteq {\cal A}_d^{>f^2(d)}$. By Lemma 2, ${\cal W}_d = {\cal
A}_d^{>f^2(d)}$ is a closed subset of ${\cal A}_d$. \hfill{$\Box$}

\medskip

\bcen{\bf 3. Morita equivalence}\ecen

\medskip

Now we study the changes of the rank of a wild algebra under
Morita equivalence and factor algebra. The following result
implies that for the proof of the Wild-Rank Conjecture it suffices
to show it for all basic algebras.

\medskip

{\bf Theorem 2.} {\it If a $d$-dimensional wild algebra $A$ is
Morita eqivalent to a basic algebra $B$ then $r_A \leq d \cdot
r_B$.}

\medskip

{\bf Proof.} Suppose $A= \oplus _{i=1}^m n_iP_i$ with $n_i \geq 1$
and $P_i, 1 \leq i \leq m,$ being the nonisomorphic indecomposable
projective $A$-modules. Let $P=\oplus _{i=1}^m P_i$. Then $B \cong
\End_A(P)^{op}$. Consider the evaluation functor $e_P=\Hom_A(P,-):
\mod A \rightarrow \mod B$. Note that $e_P$ is an equivalence of
categories with quasi-inverse $P \otimes _B -$ (cf. [2; Corollary
II.2.6.] and [1; Theorem 22.2]). Since $B$ is wild, there is a
$B$-$k \langle x,y \rangle$-bimodule $M$ which is free of rank
$r_B$ over $k \langle x,y \rangle$ such that the functor $M
\otimes _ {k \langle x,y \rangle} -$ from $\mod k \langle x,y
\rangle$ to $\mod B$ preserves indecomposability and isomorphism
classes. Note that $P$ is also projective over $B$. Decompose $P$
as the direct sum of the indecomposable projective right
$B$-modules, set $P= \oplus _{i=1}^t Q_i$. For $Q_i$ there is a
projective right $B$-module $Q_i'$ such that $Q_i \oplus Q_i' =B$.
Thus there is a projective right $B$-module $P'$ such that $P
\oplus P' =B^t$. Further $(P \otimes _B M) \oplus (P' \otimes _B
M)=B^t \otimes _B M$ which is free of rank $t \cdot r_B \leq
\dim_kP \cdot r_B \leq \dim_kA \cdot r_B= d \cdot r_B$. Since $P
\otimes _B M$ is finitely generated projective over $k \langle x,y
\rangle$, by [5; Theorem 1.4.1], it is free over $k \langle x,y
\rangle$. Moreover, its rank is at most $d \cdot r_B$. Consider
the composition $P \otimes _B M \otimes _ {k \langle x,y \rangle}
-$, we have $r_A \leq d \cdot r_B$. \hfill{$\Box$}

\medskip

From now on, unless stated otherwise, {\it we assume that all
algebras are basic.} Thus any algebra $A$ can be written as $kQ/I$
where $Q$ is the Gabriel quiver of $A$ and $I$ is an admissible
ideal of the path algebra $kQ$. For a quiver $Q$ we denote by
$Q_0$ (resp. $Q_1$) the set of vertices (resp. arrows) of $Q$. The
next result implies that for the proof of the Wild-Rank Conjecture
it suffices to show it for all minimal wild algebras. Here {\it
minimal wild} means no proper factor algebra is wild.

\medskip

{\bf Lemma 3.} {\it If $I$ is an ideal of an algebra $A$ and $A/I$
is wild then $r_A \leq r_{A/I}$.}

\medskip

{\bf Proof.} If $M$ is a finitely generated $A/I$-$k \langle x,y
\rangle$-bimodule which is free of rank $r_{A/I}$ over $k \langle
x,y \rangle$ such that the functor $M \otimes _ {k \langle x,y
\rangle} -$ from $\mod k \langle x,y \rangle$ to $\mod A/I$
preserves indecomposability and isomorphism classes, then $M$ is
also a finitely generated $A$-$k \langle x,y \rangle$-bimodule
which is free of rank $r_{A/I}$ over $k \langle x,y \rangle$ such
that the functor $M \otimes _ {k \langle x,y \rangle} -$ from
$\mod k \langle x,y \rangle$ to $\mod A$ preserves
indecomposability and isomorphism classes. \hfill{$\Box$}

\medskip

\bcen{\bf 4. Covering criterion}\ecen

\medskip

In this section, we shall provide a covering criterion which can
be effectively applied to provide an anticipated upper bound for
the rank of a concrete wild algebra. For the knowledge of Galois
covering theory we refer to [4, 12, 27].

A {\it minimal wild concealed algebra} means a concealed algebra
of a minimal wild hereditary algebra. Unless stated otherwise, the
{\it minimal} in {\it minimal wild hereditary algebra} or {\it
minimal wild concealed algebra} is always in the sense of [21].
First of all, we provide upper bounds for the ranks of some
strictly wild subcategories in the module categories over minimal
wild concealed algebras.

\medskip

{\bf Lemma 4.} {\it The ranks of all minimal wild hereditary
algebras are bounded by a fixed number.}

\medskip

{\bf Proof.} Note that the underlying diagrams of the quivers of
all minimal wild hereditary algebras are listed in [21; p.443].
Denote by $|Q|$ the underlying diagram of the quiver $Q$. Then
there are at most $2^{|Q_1|}$ quivers with underlying diagram
$|Q|$. Thus (up to isomorphism) there are finitely many minimal
wild hereditary algebras. \hfill{$\Box$}

\medskip

Let $A=kQ/I$. For an $A$-module $M$ we define its {\it support}
$\Supp(M)$ to be the subset of $Q_0$ consisting of those $x \in
Q_0$ satisfying $M(x) \not= 0$. An $A$-module $M$ is called {\it
sincere} if $\Supp(M)=Q_0$.

\medskip

{\bf Lemma 5.} {\it The ranks of all minimal wild concealed
algebras are bounded by a fixed number.}

\medskip

{\bf Proof.} It is enough to show that (up to isomorphism) there
are only finitely many minimal wild concealed algebras. This is
clear by [32; 33]. Here we give some details. Let $A$ be a minimal
wild concealed algebra of type $H$. Let $T= \oplus ^n_{i=1}T_i$ be
a preprojective tilting $H$-module such that $A=\End_H(T)$. Then
$T_i=\tau^{-m_i}P_i$ for some indecomposable projective $H$-module
$P_i$ and some nonnegative integer $m_i$. Here $\tau$ denotes
Auslander-Reiten translation. Thus $T=\tau^{-min\{m_i|1 \leq i
\leq n\}}T_1$ with $T_1=P \oplus \tau^{-1} T_2$, where $P$ is a
projective $H$-module and $\tau^{-1} T_2$ has no projective direct
summand. By [31; p.76, (6)]) we have $\Ext^1_H(T_1,T_1)=0$. Thus
$T_1$ is still a preprojective tilting $H$-module. By [2;
Proposition 1.9 (b)] we have $\End_H(T_1)=\End_H(T)=A$. Let $P=He$
and $H'=H/\langle e \rangle$ where $\langle e \rangle$ is the
two-sided ideal of $H$ generated by $e$. Then $\Hom_H(P,T_2) =
\Hom_H(P,\tau\tau^{-1}T_2) = D\Ext^1_H(\tau^{-1}T_2,P) =0$. Thus
$T_2$ is an $H'$-module. In particular $T_2$ is a non-sincere
preprojective $H$-module. Since there are only finitely many
non-sincere indecomposable preprojective $H$-modules (cf. [23;
Corollary 3.9]), there are only finitely many square-free
preprojective tilting $H$-modules with projective summands.
Therefore there are only finitely many minimal wild concealed
algebras of type $H$. By the proof of Lemma 4 the number of
minimal wild hereditary algebras is finite, so is the number of
minimal wild concealed algebras. \hfill{$\Box$}

\medskip

Denote by $(\mod A)_s$ the full subcategory of $\mod A$ consisting
of all $A$-modules whose indecomposable direct summands are all
sincere. Note that this notation is different from that in [10,
16].

\medskip

{\bf Lemma 6.} {\it If $A=kQ/I$ is a strictly wild algebra and
$A/\langle e_i \rangle$ is not strictly wild for any primitive
idempotent corresponding to the vertex $i$ in $Q_0$, then
$(\mbox{\rm mod}A)_s$ is strictly wild.}

\medskip

{\bf Proof.} The proof is almost the same as that of [16; Lemma
(3.1)]. Denote by $\mathbb{K}_3$ the quiver with two vertices
$1,2$ and tree arrows $\az, \bz, \gz$. First of all, there is a
fully faithful exact functor ${\cal F}: \mod k{\mathbb{K}}_3 \lra
\mod k \la x,y \ra$, which is defined by sending
$(V_1,V_2;\az,\bz,\gz)$ to

\noindent $((V_1 \oplus V_2)^7;
{\tiny \left[\begin{array}{ccccccc} 0&1&0&0&0&0&0\\0&0&1&0&0&0&0\\0&0&0&1&0&0&0\\
0&0&0&0&1&0&0\\0&0&0&0&0&1&0\\0&0&0&0&0&0&1\\0&0&0&0&0&0&0\end{array}
\right]},
{\tiny \left[\begin{array}{ccccccc} 0&0&0&0&0&0&0\\1&0&0&0&0&0&0\\
\sz&1&0&0&0&0&0\\
0&\delta&1&0&0&0&0\\0&0&\az'&1&0&0&0\\0&0&0&\bz'&1&0&0\\0&0&0&0&\gz'&1&0
\end{array} \right]})$ where the entries of two matrices all are $2 \times 2$
matrices and $\sz={\tiny
\left[\begin{array}{cc}1&0\\0&0\end{array} \right]}, \delta={\tiny
\left[\begin{array}{cc}0&0\\0&1\end{array} \right]}, \az'={\tiny
\left[\begin{array}{cc}0&0\\ \az&0\end{array} \right]},
\bz'={\tiny \left[\begin{array}{cc}0&0\\ \bz&0\end{array}
\right]}$ and $\gz'={\tiny \left[\begin{array}{cc}0&0\\
\gz&0\end{array} \right]}$. Moreover, there is also a fully
faithful exact functor ${\cal G}: \mod k \la x,y \ra \lra \mod
k{\mathbb{K}}_3$ which is defined by sending $(V;x,y)$ to
$(V,V;1,x,y).$ Since $A$ is strictly wild, there exists a fully
faithful exact functor ${\cal H}: \mod k{\mathbb{K}}_3 \lra \mod
A$. By assumption, we know that Supp$({\cal H}(S_1)) \cup
$Supp$({\cal H}(S_2)) = Q_0$, where $S_i$ is the simple
$k{\mathbb{K}}_3$-module corresponding to vertex $i$. It is easy
to see that both ${\cal G}{\cal F}(S_1)$ and ${\cal G}{\cal
F}(S_2)$ are sincere $k{\mathbb{K}}_3$-modules, i.e. for each $i$,
${\cal G}{\cal F}(S_i)$ is an extension of ${S_1}^{m_i}$ by
${S_2}^{n_i}$ for some positive integers $m_i$ and $n_i$. Hence
${\cal H}{\cal G}{\cal F}(S_1)$ and ${\cal H}{\cal G}{\cal
F}(S_2)$ are sincere $A$-modules. Since the functor ${\cal H}{\cal
G}{\cal F}$ is fully faithful exact, it preserves
indecomposability. Hence each indecomposable direct summand of
each $A$-module in $\Im {\cal H}{\cal G}{\cal F}$ is an image of a
module in $\mod k{\mathbb{K}}_3$. Thus all $A$-modules in $\Im
{\cal H}{\cal G}{\cal F}$ are contained in $(\mbox{\rm mod}A)_s$.
Finally ${\cal H}{\cal G}{\cal F}{\cal G}$ defines a strictly wild
functor from $\mod k \la x,y \ra$ to $(\mbox{\rm mod}A)_s$.
\hfill{$\Box$}

\medskip

The constant $b$ in the next lemma is very important, and it will
appear frequently.

\medskip

{\bf Lemma 7.} {\it The ranks of} $(\mod A)_s$ {\it where $A$ runs
through all minimal wild concealed algebras are bounded by a fixed
number. Suppose $b$ is the smallest bound.}

\medskip

{\bf Remark 2.} It should be interesting to evaluate the number
$b$.

\medskip

{\bf Proof.} It follows from [22; Corollary 2.2] that $\mod A$ is
strictly wild. It is well-known that the   minimal wild concealed
algebras are minimal wild in the sense of [21] (cf. [33; p.146]).
By Lemma 6, we know $(\mod A)_s$ is strictly wild as well. By the
proof of Lemma 5, we know there are only finitely many minimal
wild concealed algebras. \hfill{$\Box$}

\medskip

A quiver with relations $(Q,I)$ is called a {\it factor quiver} of
a quiver with relations $(Q',I')$ if $Q_0$ is a subset of $Q'_0$,
$Q_1$ is a subset of the subset of $Q'_1$ obtained from $Q'_1$ by
excluding all the arrows starting or ending at some vertex in
$Q'_0 \backslash Q_0$, and $I$ is the admissible ideal of $kQ$
obtained from $I'$ by replacing each arrow in $Q'_1 \backslash
Q_1$ in each element of $I'$ by zero (cf. [16]). Note that in this
case $kQ/I$ is a factor algebra of $kQ'/I'$. A Galois covering of
quiver with relation $\pi : (Q',I') \rightarrow (Q,I)$ is said to
be {\it wild concealed} if there is a finite factor quiver
$(\tilde{Q},\tilde{I})$ of $(Q',I')$ such that
$k\tilde{Q}/\tilde{I}$ is a minimal wild concealed algebra. The
following result including its proof is a modification of [10;
Proposition I.10.6].

\medskip

{\bf Lemma 8.} {\it Let $\pi : (Q',I') \rightarrow (Q,I)$ be a
Galois covering of quiver with relations with torsion-free Galois
group $G$ and $(\tilde{Q},\tilde{I})$ a finite factor quiver of
$(Q',I')$. Then

(1) The restriction} $F_{\lambda} : (\mod k\tilde{Q} /
\tilde{I})_s \rightarrow \mod kQ / I$ {\it preserves
indecomposability and isomorphism classes.

(2) There is a finitely generated $kQ/I$-$k\tilde{Q} /
\tilde{I}$-bimodule $M$ which is free of rank $|\tilde{Q}_0|$ over
$k\tilde{Q} / \tilde{I}$ such that on} $(\mod k\tilde{Q} /
\tilde{I})_s$, $F_{\lambda} \cong M \otimes _{k\tilde{Q} /
\tilde{I}}-$.

\medskip

{\bf Proof.} (1) $F_{\lambda}$ preserves indecomposability:
Suppose $N$ is an indecomposable in $(\mod k\tilde{Q} /
\tilde{I})_s$. Then we consider $N$ as a $kQ'/I'$-module. By [12;
Lemma 3.5], it suffices to show that $^gN \ncong N$ for $1 \neq g
\in G$. If $1 \neq g$ then, since $G$ is torsion-free, $(^g
\tilde{Q})_0 \neq \tilde{Q}_0$. Hence $\Supp(^g N) \neq \Supp(N)$.
Thus $^g N \ncong N$.

$F_{\lambda}$ preserves isomorphism classes: Let $F_{\lambda}(N_1)
\cong F_{\lambda}(N_2)$. Let $N_j= \oplus ^{n_j}_{i=1}N_{ji}$ be
the direct sum decomposition of $N_j \in (\mod k\tilde{Q} /
\tilde{I})_s$, $j=1,2$, into indecomposables. Then, by the
paragraph above and Krull-Schmidt theorem, we have $n_1=n_2$ and
$F_{\lambda}(N_{1i}) \cong F_{\lambda}(N_{2t_i}), 1 \leq t_i \leq
n_1, i=1,...,n_1$. Considering $N_{ji}, j=1,2, i=1,...,n_1$ as
$kQ'/I'$-module. By [12; Lemma 3.5], we have $N_{1i} \cong
^{g_i}\!\!N_{2t_i}$ for some $g_i \in G$ and $i=1,...,n_1$. Thus
$\tilde{Q}_0 = \Supp(N_{1i}) = \Supp(^{g_i}\!N_{2t_i}) =
^{g_i}\!\!\tilde{Q}_0$. Since $G$ is torsion-free, we have $g_i
=1$ and $N_{1i} \cong N_{1t_i}, i=1,...,n_1$. Hence $N_1 \cong
N_2$.

(2) The $kQ/I$-$k\tilde{Q} / \tilde{I}$-bimodule $M$: Define $M$
to be the free $k\tilde{Q} / \tilde{I}$-module $\oplus_{i \in
\tilde{Q}_0}b_i(k\tilde{Q} / \tilde{I})$ with free basis $\{b_i |
i \in \tilde{Q}_0\}$. We define a left $kQ/I$-module structure on
$M$ as follows: Let $i \in Q_0$, $s \in \tilde{Q}_0$ and $\sz \in
k\tilde{Q} / \tilde{I}$. We denote by $e_s$ the idempotent of
$k\tilde{Q}$ corresponding to $s$, and we set $e_i(b_s\sz)=\left\{
\begin{array}{ll} b_s(e_s\sz) & \mbox{if $\pi(s)=i$,}\\ $0$ &
\mbox{otherwise.} \end{array} \right.$ Suppose $\az : i
\rightarrow j$ is an arrow in $Q$. If $s \in \tilde{Q}_0$ with
$\pi(s)=i$ and $\tilde{\az} : s \rightarrow t$ is an arrow in
$\tilde{Q}$ with $\pi(s)=i$ and $\pi(\tilde{\az})=\az$ then we
define $\az(b_s\sz)=b_t(\tilde{\az}\sz)$, and set $\az(b_s\sz)=0$
otherwise. We claim that this is a $kQ/I$-module action: Suppose
$\rho \in I$. Note that every relation is the sum of minimal and
zero relations (cf. [27]). For the proof of $\rho(b_s\sz)=0$ for
$\sz \in k\tilde{Q} / \tilde{I}$ it suffices to show it for
minimal or zero relation $\rho \in I$. We assume $\rho \in
e_j(kQ)e_i$ for $i,j \in Q_0$. If there is no $s \in \tilde{Q}_0$
such that $\pi(s)=i$ then we have $\rho(b_s\sz)=0$. If there is $s
\in \tilde{Q}_0$ such that $\pi(s)=i$ then there is $\rho ' \in I'
\cap e_t(kQ')e_s$ such that $\pi(\rho ')=\rho$. By replacing each
arrow in $Q'_1 \backslash \tilde{Q}_1$ by zero we obtain
$\tilde{\rho} \in \tilde{I} \cap e_t(k \tilde{Q})e_s$ from $\rho
'$. Clearly $\rho(b_s\sz)=b_t(\tilde{\rho}\sz)=0$.

Now let $N \in \mod k\tilde{Q} / \tilde{I}$, we will show that
$F_{\lambda}(N)=M \otimes _{k\tilde{Q} / \tilde{I}} N$
canonically. Since for any arrow $\tilde{\az} \in \tilde{Q}$ we
have that $(b_s\tilde{\az}) \otimes N = b_s \otimes (\tilde{\az}
N) \subseteq b_s \otimes N$, the module $M \otimes _{k\tilde{Q} /
\tilde{I}} N$ has underlying space $\oplus _{s \in
\tilde{Q}_0}(b_s \otimes N)$. Let $i \in Q_0$. If $\pi(s) \neq i$
then $e_i(b_s \otimes N)=0$. If $\pi(s) = i$ then $e_i(b_s \otimes
N)= (b_se_s) \otimes N = b_s \otimes e_sN = b_s \otimes N(s)$. So
we may identify $e_i(M \otimes N)$ with $(F_{\lambda}(N))(i)=
\oplus _{\pi(s) = i}N(s)$. Now consider the action of an arrow
$\az : i \rightarrow j$ in $Q$. Let $\tilde{\az} : s \rightarrow
t$ be an arrow in $\tilde{Q}$ with $\pi(s)=i$,
$\pi(\tilde{\az})=\az$ and hence $\pi(t)=j$. Then $\az(b_s \otimes
N) = (b_t \tilde{\az}) \otimes N = b_t \otimes (\tilde{\az} N) =
b_t \otimes (\tilde{\az} e_sN) = b_t \otimes (\tilde{\az} N(s)) =
b_t \otimes N({\tilde{\az}})( N(s))$ and this is just the action
of $\az$ on the space $(F_{\lambda}(N))(i).$ \hfill{$\Box$}

\medskip

{\bf Theorem 3.} (covering criterion) {\it Let $A=kQ/I$ be a wild
algebra and $\pi: (Q',I') \rightarrow (Q,I)$ a wild concealed
Galois covering of quivers with relations with torsion-free Galois
group. Then $r_A \leq 10b$.}

\medskip

{\bf Proof.} Let $(\tilde{Q}, \tilde{I})$ be a finite factor
quiver of $(Q',I')$ such that $k \tilde{Q} / \tilde{I}$ is a
minimal wild concealed algebra. By Lemma 7, there is a finitely
generated $k\tilde{Q}/\tilde{I}$-$k \langle x,y \rangle$-bimodule
$M_1$ which is free of rank at most $b$ over $k \langle x,y
\rangle$ such that the functor $M_1 \otimes _ {k \langle x,y
\rangle} -$ from $\mod k \langle x,y \rangle$ to $(\mod
k\tilde{Q}/\tilde{I})_s$ preserves indecomposability and
isomorphism classes. By Lemma 8, there is a finitely generated
$kQ/I$-$k\tilde{Q}/\tilde{I}$-bimodule $M_2$ which is free of rank
$|\tilde{Q}_0|$ over $kQ/I$ such that on $\mod(\tilde{Q},
\tilde{I})_s$ the pushdown functor $F_{\lambda} \cong M_2 \otimes
_{k\tilde{Q}/\tilde{I}} -$ preserves indecomposability and
isomorphism classes. Consider the composition $M_2 \otimes
_{k\tilde{Q}/\tilde{I}}M_1 \otimes _ {k \langle x,y \rangle} -$,
we have $r_A \leq \rank (M_2 \otimes M_1) \leq |\tilde{Q}_0| \cdot
b \leq 10b.$ \hfill{$\Box$}

\medskip

According to Theorem 2 and 3, we reformulate the Wild-Rank
Conjecture as follows:

\medskip

{\bf Wild-Rank Conjecture.} {\it Let $A$ be a $d$-dimensional
(unnecessarily basic) wild algebra. Then $r_A \leq 10bd$.}

\medskip

{\bf Basic-Wild-Rank Conjecture.} {\it Let $A$ be a
$d$-dimensional basic wild algebra. Then $r_A \leq 10b$.}

\medskip

Clearly, Basic-Wild-Rank Conjecture $\Rightarrow$ Wild-Rank
Conjecture $\Rightarrow$ Tame-Open Conjecture.

\medskip

\bcen{\bf 5. Applications of the covering criterion}\ecen

\medskip

How to support the Basic-Wild-Rank Conjecture? For concrete
algebras, our covering criterion is very effective. Indeed, for a
concrete basic wild algebra $A$ given by quiver with relations
$(Q,I)$, we can find a minimal wild factor algebra $B$ of $A$.
Usually either $B$ is itself a minimal wild concealed algebra or
there is an algebra $C \cong B$ such that $C$ admits a wild
concealed Galois covering with torsion-free Galois group. Thus we
can apply the covering criterion to the algebra $C$.

By the covering criterion, we know the Basic-Wild-Rank Conjecture
holds for all well-known wild algebras such as wild local
algebras, wild two-point algebras, wild radical square zero
algebras, wild finite $p$-group algebras, wild three-point
algebras whose quiver is system quiver (cf. [30, 17, 15, 16, 26]).
This implies that all three conjectures are much reliable.

Certainly one can list many propositions analogous to the
following one.

\medskip

{\bf Proposition.} {\it Let $A$ be a $d$-dimensional wild local
algebra (resp. wild two-point algebra, wild radical square zero
algebra). Then $r_A \leq 10b$.}

\medskip

{\bf Proof.}  Up to duality and isomorphism, $A$ has a minimal
wild factor algebra $B$ appearing in the list of [30; p.283]
(resp. [17; Table W], [15; p.98] or [16; p.290]). Check case by
case we know that either $B$ is itself a minimal wild concealed
algebra or there is an algebra $C \cong B$ such that $C$ admits a
wild concealed Galois covering with torsion-free Galois group.
\hfill{$\Box$}

\medskip

\bcen{\bf ACKNOWLEDGEMENT}\ecen

\medskip

The author is grateful to Otto Kerner for his explanations of some
results in the representation theory of wild tilted algebras.

\end{document}